\documentclass[10pt]{amsart}
\usepackage{amssymb,latexsym}
\usepackage{graphics,psfig}
\newcommand{\psdraw}[4]{\begin{array}{c} \hspace{-1mm}
\raisebox{-4pt}{\psfig{figure=#1.ps,width=#2,height=#3,angle=#4}}
\end{array}}

\numberwithin{equation}{section}
\newtheorem{theorem}{Theorem}[section]
\newtheorem{proposition}[theorem]{Proposition}
\newtheorem{corollary}[theorem]{Corollary}

\newtheorem{lemma}[theorem]{Lemma}

\begin{document}

\pagenumbering{arabic}
\pagestyle{headings}
\def\sof{\hfill\rule{2mm}{2mm}}
\def\ls{\leq}
\def\gs{\geq}
\def\SS{\mathcal S}
\def\qq{{\bold q}}
\def\txx{{\frac1{2\sqrt{x}}}}
\def\II{\mathcal{I}}
\def\figa{$\psdraw{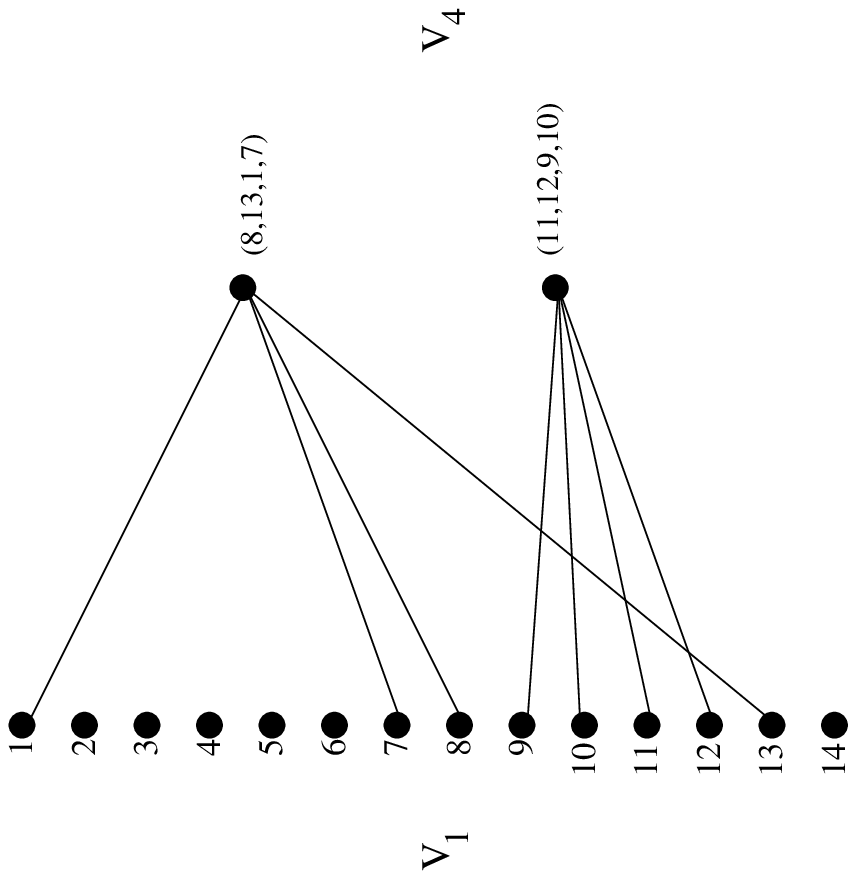}{80mm}{70mm}{270}$}
\def\figb{$\psdraw{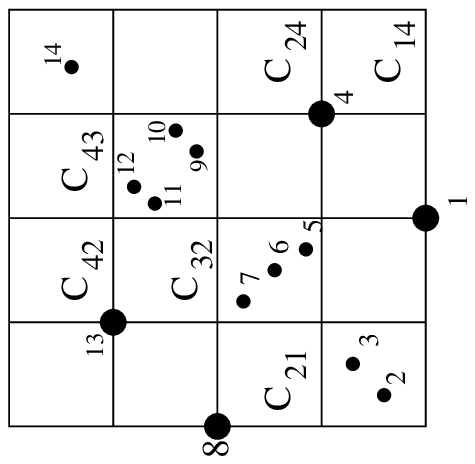}{37mm}{30mm}{270}$}
\def\figc{$\psdraw{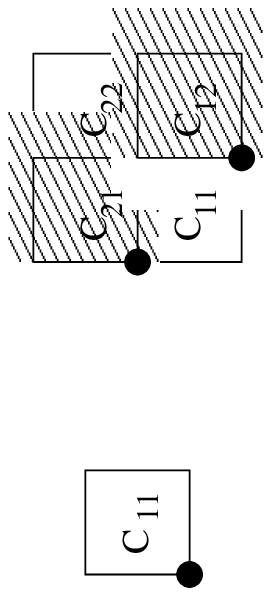}{30mm}{15mm}{270}$}
\def\figd{$\psdraw{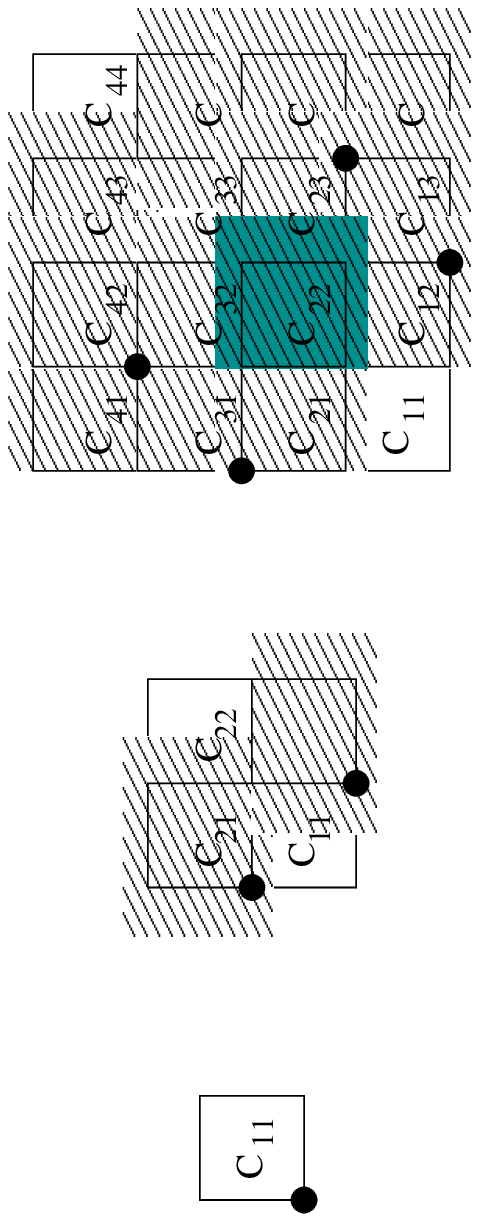}{70mm}{25mm}{270}$}

\title{Counting occurrences of $3412$ in an involution}
\maketitle
\begin{center}
Toufik Mansour
\end{center}

\begin{center}
{Department of Mathematics, Haifa University, 31905 Haifa, Israel

        {\tt toufik@math.haifa.ac.il} }
\end{center}
%
\section*{Abstract}
We study the generating function for the number of involutions on
$n$ letters containing exactly $r\gs0$ occurrences of $3412$. It
is shown that finding this function for a given $r$ amounts to a
routine check of all involutions on $2r+1$ letters.

\noindent {\sc 2000 Mathematics Subject Classification}: Primary
05A05, 05A15; Secondary 05C90

\section{Introduction}
{\bf Permutations}. Let $\pi\in S_n$ and $\tau\in S_{m}$ be two
permutations. An {\it occurrence\/} of $\tau$ in $\pi$ is a
subsequence $1\ls i_1<i_2<\dots<i_{m}\ls n$ such that
$(\pi(i_1),\dots,\pi(i_{m}))$ is order-isomorphic to $\tau$; in
such a context, $\tau$ is usually called a {\it pattern}.

In the last decade much attention has been paid to the problem of
finding the numbers $s_{\sigma}^r(n)$ for a fixed $r\geq 0$ and a
given pattern $\tau$ (see
\cite{AlFr00,At99,Bo97b,Bo97a,Bo98,ChWe99,Ma00,MaVa01,NoZe96,
Ro99,SiSc85,St94,St96,We95}). Most of the authors consider only
the case $r=0$, thus studying permutations {\it avoiding\/} a
given pattern. Only a few papers consider the case $r>0$, usually
restricting themselves to patterns of length $3$.  Using two
simple involutions (\emph{reverse} and \emph{complement}) on
$\SS_n$ it is immediate that with respect to being
equidistributed, the six patterns of length three fall into the
two classes $\{123,321\}$ and $\{132,213,231,312\}$. Noonan
\cite{No96} proved that $s_{123}^1(n)=\frac 3n\binom{2n}{n-3}$. A
general approach to the problem was suggested by Noonan and
Zeilberger \cite{NoZe96}; they gave another proof of Noonan's
result, and conjectured that
        $s_{123}^2(n)=\frac{59n^2+117n+100}{2n(2n-1)(n+5)}\binom{2n}{n-4}$
and $s_{132}^1(n)=\binom{2n-3}{n-3}$. The first conjecture was
proved by Fulmek \cite{Fulmek} and the latter conjecture was
proved by B\'ona in \cite{Bo98}. A conjecture of Noonan and
Zeilberger states that $s_{\sigma}^r(n)$ is $P$-recursive in $n$
for any $r$ and $\tau$. It was proved by B\'ona \cite{Bo97c} for
$\sigma=132$. Mansour and Vainshtein \cite{MaVa01} suggested a new
approach to this problem in the case $\sigma=132$, which allows
one to get an explicit expression for $s_{132}^r(n)$ for any given
$r$. More precisely, they presented an algorithm that computes the
generating function $\sum_{n\geq0} s_{132}^r(n)x^n$ for any
$r\geq0$. To get the result for a given $r$, the algorithm
performs certain routine checks for each element of the symmetric
group $\SS_{2r}$. The algorithm has been implemented in C, and
yields explicit results for $1\leq r\leq 6$.

{\bf Involutions}. An {\em involution} $\pi$ is a permutation such
that $\pi=\pi^{-1}$; let $\II_n$ denote the set of all
$n$-involutions. An {\em even} (resp. {\em odd}) involution $\pi$
is an involution such that the number of occurrences of the
pattern $21$ in $\pi$ is an even (resp. odd) number.

Several authors have given enumerations of sets of involutions
which avoid certain patterns. In \cite{Regev} Regev provided an
asymptotic formula for $\II_n(12\ldots k)$ and showed that
$\II_n(1234)$ is enumerated by the $n$th Motzkin number
$M_n=\sum_{i=0}^{[n/2]}\binom{n}{2i}C_i$. In \cite{Gessel} Gessel
enumerated $\II_n(12\ldots k)$. In \cite{GouyouBYoung}
Gouyou--Beauchamps gave an entirely bijective proof of some very
nice exact formulas for $\II_n(12345)$ and $\II_n(123456)$.

Pattern-avoiding involutions have also been related with other
combinatorial objects. Gire \cite{GireThese} has established a
one-to-one correspondence between 1-2 trees having $n$ edges and
$\SS_n(321,3\overline{1}42)$ (this is the set of $n$-permutations
avoiding patterns $321$ and $231$, except that the latter is
allowed when it is a subsequence of the pattern $3142$). On the
other hand, Guibert \cite{GuibertThese} has established bijections
between 1-2 trees having $n$ edges and each of the sets
$\SS_n(231,4\overline{1}32)$, $\II_n(3412)$, and $\II_n(4321)$
(and therefore with $\II_n$(1234), by transposing the
corresponding Young tableaux obtained by applying the
Robinson-Schensted algorithm). Also, Guibert \cite{GuibertThese}
has established a bijection between $\II_n(2143)$ and
$\II_n(1243)$. More recently, Guibert, Pergola, and Pinzani
\cite{GPP} have given a one-to-one correspondence between 1-2
trees having $n$ edges and $\II_n(2143)$. It follows that all
these sets are enumerated by the $n$th Motzkin number $M_n$. It
remains an open problem to prove the conjecture of Guibert (in
\cite{GuibertThese}) that $\II_n(1432)$ is also enumerated by the
$n$th Motzkin number $M_n$. This conjecture proved by~
Jaggard~\cite{Jag}. Recently, Egge~\cite{Egge} gave enumerations
and generating functions for the number of involutions in
$\II_n(3412)$ which avoid various sets of additional patterns.
Egge and Mansour~\cite{EM} refined Egge's results by studying
generating functions for the number of even and odd involutions in
$\II_n(3412)$, and they studied the generating functions for the
number of even and odd involutions in $\II_n$ containing the
pattern $3412$ exactly once. For example, they proved the
following theorem.

\begin{theorem}{\rm(see \cite[Proposition 7.7]{EM})}
The generating function for the number of involutions in $\II_n$
which contain the pattern $3412$ exactly once is given by
$$\frac{2x-1}{2x^2(1-x)}+\frac{1-2x-2x^2}{2x^2\sqrt{1-2x-3x^2}}.$$
\end{theorem}

\begin{center}
\begin{table}[h]
\begin{tabular}{|c|c|c|c|c|c|c|c|c|c|c|c|c|c|}\hline
  $r\backslash n$ & $0$ & $1$ & $2$ & $3$ & $4$ & $5$ & $6$ & $7$ & $8$ & $9$ & $10$ & $11$ &
  $12$\\\hline\hline
  $0$ & $1$ & $1$ & $2$ & $4$ & $9$ & $21$ & $51$ & $127$ & $323$ & $835$ & $2188$ & $5798$ & $15511$\\ \hline
  $1$ & $0$ & $0$ & $0$ & $0$ & $1$ & $5$  & $20$ & $70$  & $231$ & $735$ & $2289$ & $7029$ & $21384$\\ \hline
  $2$ & $0$ & $0$ & $0$ & $0$ & $0$ & $0$  & $1$ & $7$  & $37$ & $165$ & $671$ & $2563$ & $9375$\\ \hline
  $3$ & $0$ & $0$ & $0$ & $0$ & $0$ & $0$  & $1$ & $4$  & $17$ & $63$ & $236$ & $877$ & $3270$\\ \hline
  $4$ & $0$ & $0$ & $0$ & $0$ & $0$ & $0$  & $2$ & $12$  & $56$ & $220$ & $803$ & $2783$ & $9364$\\ \hline
  $5$ & $0$ & $0$ & $0$ & $0$ & $0$ & $0$  & $0$ & $2$  & $14$ & $80$ & $383$ & $1658$ & $6690$\\ \hline
  $6$ & $0$ & $0$ & $0$ & $0$ & $0$ & $0$  & $0$ & $2$  & $11$ & $51$ & $212$ & $856$ & $3402$\\ \hline
\end{tabular}
\caption{First values for the sequence $\II_{3412}^r(n)$ where
$0\leq n\leq 12$ and $0\leq r\leq 6$}\label{tab1}\vspace*{-25pt}
\end{table}
\end{center}
In this paper, we suggest a new approach to study the number of
involutions which contain the pattern $3412$ exactly $r$ times,
namely $\II_{3412}^r(n)$, which allows one to get an explicit
expression for $\II_{3412}^r(n)$ for any given $r$ (see
Table~\ref{tab1}). More precisely, we present an algorithm that
computes the generating function $\II_r(x)=\sum_{n\geq0}
\II_{3412}^r(n)x^n$ for any $r\geq0$. To get the result for a
given $r$, the algorithm performs certain routine checks for each
element of the involutions $\II_{2r+1}$. The algorithm has been
implemented in C, and yields explicit results for $0\leq r\leq 6$.
\section{Recall definitions and preliminary results}

\subsection{The case $r\geq1$}
Fix $r\geq1$. To any $\pi\in\II_n$ we assign a bipartite graph
$G_\pi$ in the following way. The vertices in one part of $G_\pi$,
denoted $V_1$, are the entries of $\pi$, and the vertices of the
second part, denoted $V_4$, are the occurrences of $3412$ in
$\pi$. Entry $i\in V_1$ is connected by an edge to occurrence
$j\in V_4$ if $i$ enters $j$. For example, let
$\pi=8\,2\,3\,13\,7\,6\,5\,1\,11\,12\,9\,10\,\,4\,14$, then $\pi$
contains $2$ occurrences of $3412$, and the graph $G_\pi$ is
presented on Figure~1.

\begin{center}
\figa

{\sc Figure 1}. Graph $G_\pi$ for
$\pi=8\,2\,3\,13\,7\,6\,5\,1\,11\,12\,9\,10\,\,4\,14$
\label{graph}
\end{center}

Let $\widetilde G$ be an arbitrary connected component of $G_\pi$,
and let $\widetilde V$ be its vertex set. We denote $\widetilde
V_1=\widetilde V\cap V_1$, $\widetilde V_4=\widetilde V\cap V_4$,
$t_1=|\widetilde V_1|$, $t_4=|\widetilde V_4|$.

\begin{lemma}\label{lem1} For any connected component $\widetilde G$ of
$G_\pi$ one has $t_1\leq 2t_4+2$.
\end{lemma}
\begin{proof}
Assume to the contrary that the above statement is not true.
Consider the smallest $n$ for which there exists $\pi\in \II_n$
such that for some connected component $\widetilde G$ of $G_\pi$
one has
        $$t_1>2t_4+2. \eqno{(\ast)}$$
Evidently, $\widetilde G$ contains more than one vertex, since
otherwise $t_1=1$ and $t_4=0$, which contradicts $(\ast)$. Let $l$
be the number of leaves in $\widetilde G$ (recall that a leaf is a
vertex of degree~$1$). Clearly, all the leaves belong to
$\widetilde V_1$; the degree of any other vertex in $\widetilde
V_1$ is at least $2$, while the degree of any vertex in
$\widetilde V_4$ equals $4$. Calculating the number of edges in
$\widetilde G$ by two different ways, we get $l+2(t_1-l)\ls 4t_4$,
which together with ($\ast$) gives $l>4$, so  there exist at least
five leaves in $\widetilde V_1$. In the case $r=1$ there only one
occurrence of length $4$, so there no five leaves in $V_1$,
therefore we can assume that $r\geq2$.

Let $a\in\widetilde V_4$ and let us consider the following cases
corresponding to the number of leaves in $V_1$ incident to $a$:
\begin{itemize}
\item[{\rm(i)}] If $a$ incident to one leave $x\in V_1$ exactly, then
the graph $\widetilde G'$ which obtained from $\widetilde G$ by
deleting $x$ and $a$ satisfies $t_1'+1=t_1$ and $t_4'+1=t_4$,
hence by ($\ast$) we have that $t_1'>2t_4'+3$, a contradiction to
the minimality of $n$.

\item[{\rm(ii)}] If $a$ incident to two leaves $x,y\in V_1$
exactly, then the graph $\widetilde G'$ which obtained from
$\widetilde G$ by deleting $x,y$ and $a$ satisfies $t_1'+2=t_1$
and $t_4'+1=t_4$, hence by ($\ast$) we have that $t_1'>2t_4'+2$, a
contradiction to the minimality of $n$.

\item [{\rm(iii)}] If there exist four leaves incident to $a$, then
the graph $\widetilde G$ does not connect component, a
contradiction.

\item[{\rm(iv)}] By Cases(i)-(iii) we can assume that every vertex
in $V_4$ incident to three leaves exactly. Consider $a,b\in V_4$
with the leaves $a^1,a^2,a^3\in V_1$ and $b^1,b^2,b^3\in V_1$ such
that there exists a vertex $u\in V_1$ incident to $a$ and $b$.
Therefore, if we consider all the involutions $\pi$ of length at
most $7$ with two occurrences $a$ and $b$ of $3412$ and the
corresponding graph $G_\pi$ is connected component, then we see
that each involution has two occurrences $x_1x_2x_3x_4$ and
$y_1y_2y_3y_4$ of $3412$ such that
$|\{x_1,x_2,x_3,x_4\}\cap\{y_1,y_2,y_3,y_4\}|=2$, a contradiction
for that every vertex in $V_4$ incident to three leaves exactly
(more precisely, in $\II_k$, $k=0,1,2,\ldots,5$, there no
involutions which contain $3412$ exactly twice; in $\II_6$ there
exists one involution which contain $3412$ exactly twice, namely
$351624$; and in $\II_7$ we have $7$ involutions which contain
$3412$ exactly twice, namely $1462735$, $3516247$, $3517264$,
$3614725$, $3617524$, $4261735$, and $4631725$. Each of these
involutions has two occurrences $x_1x_2x_3x_4$ and $y_1y_2y_3y_4$
of $3412$ such that there
$|\{x_1,x_2,x_3,x_4\}\cap\{y_1,y_2,y_3,y_4\}|=2$).
\end{itemize}

Hence, using Cases (i)-(iv) we get the desired result.
\end{proof}

Denote by $G_\pi^1$ the connected component of $G_\pi$ containing
entry $1$. Let $\pi(i_1),\dots,\pi(i_s)$ be the entries of $\pi$
belonging to $G_\pi^1$, and let $\sigma=\sigma_\pi\in \II_s$ be
the corresponding involution ($\sigma_\pi$ is an involution for
any involution $\pi$) . We say that $\pi(i_1),\dots,\pi(i_s)$ is
the {\it kernel\/} of $\pi$ and denote it $\ker\pi$; $\sigma$ is
called the {\it shape\/} of the kernel, or the {\it kernel shape},
$s$ is called the {\it size\/} of the kernel, and the number of
occurrences of $3412$ in $\ker\pi$ is called the {\it capacity\/}
of the kernel. For example, for
$\pi=8\,2\,3\,13\,7\,6\,5\,1\,11\,12\,9\,10\,\,4\,14$ as above,
the kernel equals $8\,13\,1\,4$, its shape is $3412$, the size
equals $4$, and the capacity equals $1$.

The following statement is implied immediately by
Lemma~\ref{lem1}.

\begin{theorem}\label{th1} Fix $r\geq1$. Let $\pi\in \II_n$ contain exactly $r$ occurrences
of $3412$, then the size of the kernel of $\pi$ is at most $2r+2$.
\end{theorem}

We say that $\rho$ is a {\it kernel involution\/} if it is the
kernel shape for some involution $\pi$. Evidently $\rho$ is a
kernel involution if and only if $\sigma_\rho=\rho$ is an
involution.

Let $\rho\in \II_s$ be an arbitrary kernel involution. We denote
by $\II(\rho)$ the set of all the involutions of all possible
sizes whose kernel shape equals $\rho$. For any $\pi\in \II(\rho)$
we define the {\it kernel cell decomposition\/} as follows. The
number of cells in the decomposition equals $s\times s$. Let
$\ker\pi=\pi(i_1),\dots,\pi(i_s)$; the {\it cell\/}
$C_{m\ell}=C_{m\ell}(\pi)$ for $1\ls\ell\ls s$ and $1\ls m\ls s$
is defined by
$$
C_{m\ell}(\pi)=\{\pi(j)|\, i_{\ell}<j<i_{\ell+1},
\;\pi(i_{\rho^{-1}(m)})<\pi(j)< \pi(i_{\rho^{-1}(m+1)})\},
$$
where $i_{s+1}=n+1$ and $\alpha_{n+1}=n+1$ for any $\alpha\in
\II_n$. If $\pi$ coincides with $\rho$ itself, then all the cells
in the decomposition are empty. An arbitrary permutation in
$\II(\rho)$ is obtained  by filling in some of the cells in the
cell decomposition. A cell $C$ is called {\it infeasible\/} if the
existence of an entry $a\in C$ would imply an occurrence of $3412$
that contains $a$, two other entries $x,y\in\ker\pi$, and another
entry $z$ of $\pi$. Clearly, all infeasible cells are empty for
any $\pi\in \II(\rho)$. All the remaining cells are called {\it
feasible\/}; a feasible cell may, or may not, be empty. The set of
all feasible cells are distributed into three subsets:
\begin{itemize}
\item[{\rm(i)}] {\it Free-cells\/}; a free-cell is a feasible cell
may, or may not, be contain an occurrence of the pattern $12$, or
contain an occurrence of the pattern $21$;

\item[{\rm(ii)}] {\it Diagonal-decreasing-cells\/}; a
diagonal-decreasing-cell is a feasible cell
$C_{ii}=\{\pi_{j_1},\ldots,\pi_{j_d}\}$ such that
$\pi_{j_1}>\cdots>\pi_{j_d}$ where $j_1<\cdots<j_d$;

\item[{\rm(iii)}] {\it Decreasing-cells\/}; a decreasing-cell is a
feasible cell $C_{ij}=\{\pi_{j_1},\ldots,\pi_{j_d}\}$ with $i\neq
j$ such that $\pi_{j_1}>\cdots>\pi_{j_d}$ where $j_1<\cdots<j_d$.
\end{itemize}

Consider the involution
$\pi=8\,2\,3\,13\,7\,6\,5\,1\,11\,12\,9\,10\,\,4\,14$. The kernel
of $\pi$ equals $8\,13\,1\,4$, its shape is $3412$. The cell
decomposition of $\pi$ contains four feasible cells:
$C_{11}=\{2,3\}$, $C_{22}=\{7,6,5\}$, $C_{33}=\{11,12,9,10\}$, and
$C_{44}=\{14\}$, see Figure~2. All the other cells are infeasible;
for example, $C_{21}$ is infeasible, since if $a\in C_{21}$, then
$a\pi'_{i_2}\pi'_{i_3}\pi'_{i_4}$ is an occurrence of $3412$ for
any $\pi'$ whose kernel is of shape $3412$. The cells $C_{11}$,
$C_{33}$, and $C_{44}$ are free-cells, and $C_{22}$ is
diagonal-decreasing-cell.

\begin{center}
\figb

{\sc Figure 2}. Kernel cell decomposition for $\pi\in \II(3412)$.
\end{center}

Given a cell $C_{ij}$ in the kernel cell decomposition, all the
kernel entries can be positioned with respect to $C_{ij}$. We say
that $x=\pi(i_k)\in\ker\pi$ lies {\it below\/} $C_{ij}$ if
$\rho(k)<i$, and {\it above\/} $C_{ij}$ if $\rho(k)\gs i$.
Similarly, $x$ lies to the {\it left\/} of $C_{ij}$ if $k<j$, and
to the {\it right\/} of $C_{ij}$ if $k \gs j$. As usual, we say
that $x$ lies to the {\it southwest\/} of $C_{ij}$ if it lies
below $C_{ij}$ and to the left of it; the other three directions,
northwest, southeast, and northeast, are defined similarly.

The following statement plays a crucial role in our
considerations.

\begin{lemma}\label{lem2} Let $\pi\in\II(\rho)$.

{\rm(i)} $C_{ij}$ is infeasible cell if and only if $C_{ji}$ is
infeasible cell.

{\rm(ii)} $C_{ij}$ is feasible cell with $i\neq j$ if and only if
$C_{ij}$ and $C_{ji}$ are decreasing-cells.

{\rm(iii)} Let $C_{ii}$ be a feasible cell; if there exits
northwest or southeast of $C_{ii}$ an occurrence of the pattern
$12$ in the kernel shape $\rho$, then $C_{ii}$ is
diagonal-decreasing-cell, otherwise $C_{ii}$ is free-cell.
\end{lemma}
\begin{proof}
(i) Holds immediately by the fact that $\pi$ is an involution.

(ii) Let $i\neq j$. Using the fact that $\pi$ is an involution we
get that $C_{ij}$ contains an occurrence of the pattern $12$ if
and only if $C_{ji}$ contains an occurrence of the pattern $12$.
Thus, if $C_{ij}$ contains an occurrence of the pattern $12$,
namely $xy$, then $C_{ji}$ contains an occurrence $uv$ of the
pattern $12$, so $xyuv$ is an occurrence of the pattern $3412$, a
contradiction. Hence, $C_{ij}$ is feasible cell if and only if
$C_{ij}$ and $C_{ji}$ are decreasing-cells.

(iii) Let $C_{ii}$ be a feasible cell such that there exists
northwest of $C_{ii}$ an occurrence of the pattern $12$ in the
kernel shape $\rho$, namely $xy$. So, if $C_{ii}$ contains an
occurrence of the pattern $12$ in $C_{ii}$, say $uv$, then $xyuv$
is an occurrence of the pattern $3412$. Hence, $C_{ii}$ is
diagonal-decreasing-cell.
\end{proof}

Similarly as the arguments in the proof of Lemma~\ref{lem2} we
have the following result.

\begin{lemma}\label{lem22}
Let $a<b$.

{\rm(i)} If $C_{ia}$ and $C_{ib}$ are decreasing-cells then every
entry of $C_{ia}$ is greater than every entry of $C_{ib}$;

{\rm(ii)} If $C_{ai}$ and $C_{bi}$ are decreasing-cells then the
entries of $C_{bi}$ are lie northwest of the entries of $C_{ib}$;

{\rm(iii)} If $C_{ab}$ is decreasing-cell, then all the cells
$C_{ij}$, with $i\neq j$, which are lie northeast of $C_{ab}$ and
$C_{ba}$ are infeasible cells.
\end{lemma}

As a consequence of Lemmas~\ref{lem2} and \ref{lem22}, we get the
following result. Let us define a partial order $\prec$ on the set
of all free-cells by saying that $C_{ml}\prec C_{m'\ell'}\ne
C_{m\ell}$ if $m\leq m'$ and $\ell\leq \ell'$. Similarly, we
define a partial order $\prec'$ and $\prec''$ on the set of all
diagonal-decreasing-cells and decreasing-cells, respectively.

\begin{lemma}\label{lem3} $\prec$, $\prec'$, and $\prec''$ are linear orders.
\end{lemma}

Lemmas~\ref{lem2}--\ref{lem3} yield immediately the following two
results.

\begin{theorem}\label{th2} Let $\widetilde G$ be a connected component of
$G_\pi$ distinct from $G_\pi^1$. Then all the vertices in
$\widetilde V_1$ belong to the same feasible cell in the kernel
cell decomposition of $\pi$.
\end{theorem}

Let $F(\rho)$ (respectively; $DD(\rho)$ and $D(\rho)$) be the set
of all free-cells (respectively; diagonal-decreasing-cells and
decreasing-cells which lie above the diagonal) in the kernel cell
decomposition corresponding to involutions in $\II(\rho)$ and let
$f(\rho)=|F(\rho)|$ (respectively; $dd(\rho)=|DD(\rho)|$ and
$d(\rho)=|D(\rho)|$). We remark that, by Lemma~\ref{lem2} and
Lemma~\ref{lem22} we get that $d(\rho)$ is a nonnegative integer
number. We denote the cells in $F(\rho)$ by
$C^1,\dots,C^{f(\rho)}$ (respectively; $DD^1,\dots,DD^{dd(rho)}$
and $D^1,\dots,D^{d(\rho)}$) in such a way that $C^i\prec C^j$
(respectively; $DD^i\prec' DD^j$ and $D^i\prec'' D^j$) whenever
$i<j$.

\begin{theorem}\label{th3} For any given sequence
$\alpha_1,\dots,\alpha_{f(\rho)}$ of arbitrary involutions, and
two sequences  $\beta_1,\dots,\beta_{dd(\rho)}$ and
$\gamma_1,\dots,\gamma_{d(\rho)}$ of arbitrary decreasing
involutions {\rm(}a {\it decreasing involution of length $n$} is
the involution $n(n-1)\ldots1${\rm)}, there exists $\pi\in
\II(\rho)$ such that the content of the free-cell $C^i$ is
order-isomorphic to $\alpha_i$, the content of the
diagonal-decreasing-cell $DD^j$ is $\beta_j$, and the content of
the decreasing-cell $D^j$ is $\gamma_j$.
\end{theorem}

\subsection{The case $r=0$}
First of all, let us describe the cell block decomposition of an
involution $\pi\in\II_n(3412)$ as follows.

\begin{proposition}{\rm(see~\cite[Proposition~2.8]{Egge})}\label{pro1}
Let $\pi\in I_n(3412)$. Then one of the following holds:

{\rm(i)} $\pi=(1,\alpha)$ where
$(\alpha_1-1,\ldots,\alpha_{n-1}-1)\in I_{n-1}(3412)$,

{\rm(ii)} there exists $t$, $2\leq t\leq n$, such that
$\pi=(t,\alpha,1,\beta)$ where
$(\alpha_1-1,\ldots,\alpha_{t-2}-1)\in I_{t-2}(3412)$ and
$(\beta_1-t,\ldots,\beta_{n-t}-t)\in I_{n-t}(3412)$.
\end{proposition}

Thus, the kernel cell decomposition of $\pi\in\II_n(3412)$ can be
defined as follows. There are two kernel shapes $\rho^1=1$ and
$\rho^2=21$. In the case $\rho^1$ there exists only one cell which
is free-cell, and in the case $\rho^2$ there exist four cells: two
are infeasible cells ($C_{21}$ and $C_{12}$), and the others are
free-cells ($C_{11}$ and $C_{22}$), see Figure~3.

\begin{center}
\figc

{\sc Figure 3}. Kernel cell decomposition for $\pi\in
\II(1)\cup\II(21)$.
\end{center}

\section{Main Theorem  and explicit results}
Let $\rho$ be a kernel involution, and let $s(\rho)$, $c(\rho)$,
$f(\rho)$, $dd(\rho)$, and $d(\rho)$ be the size of $\rho$, the
capacity of $\rho$, and the number of free-cells,
diagonal-decreasing-cells, decreasing-cells in the cell
decomposition associated with $\rho$, respectively. Denote by $K$
the set of all kernel involutions, and by $K_t$ the set of all
kernel shapes for involutions in $\II_t$. The main result of this
note can be formulated as follows.

\begin{theorem}\label{th4} For any $r\geq1$,
\begin{equation}
\II_r(x)=\sum_{\rho\in
K_{2r+2}}\left(\frac{x^{s(\rho)}}{(1-x^2)^{d(\rho)}(1-x)^{dd(\rho)}}
\sum_{r_1+\dots+r_{f(\rho)}=r-c(\rho)}
\;\prod_{j=1}^{f(\rho)}\II_{r_j}(x)\right),\label{mth}
\end{equation}
where $r_j\geq0$ for $1\leq j\leq f(\rho)$.
\end{theorem}
\begin{proof} For any $\rho\in K$, denote by $\II_r^\rho(x)$ the
generating function for the number of involutions in
$\pi\in\II_n\cap \II(\rho)$ containing exactly $r$ occurrences of
$3412$. Evidently, $\II_r(x)=\sum_{\rho\in K} \II_r^\rho(x)$. To
find $\II_r^\rho(x)$, recall that the kernel of any $\pi$ as above
contains exactly $c(\rho)$ occurrences of $3412$. The remaining
$r-c(\rho)$ occurrences of $3412$ are distributed between the
free-cells of the kernel cell decomposition of $\pi$. By
Theorem~\ref{th2}, each occurrence of $3412$ belongs entirely to
one free-cell. Besides, it follows from Theorem~\ref{th3}, that
occurrences of $3412$ in different free-cells do not influence one
another. Also, by Lemmas~\ref{lem2} and \ref{lem22} the
diagonal-decreasing-cells do not influence one another, and for
decreasing cell $C_{ij}$ we have that $|C_{ij}|=|C_{ji}|$ where
$C_{ij}$ contains arbitrary decreasing involution. Therefore,
$$
\II_r^\rho(x)=\frac{x^{s(\rho)}}{(1-x^2)^{d(\rho)}(1-x)^{dd(\rho)}}
\sum_{r_1+\dots+r_{f(\rho)}=r-c(\rho)}\;\prod_{j=1}^{f(\rho)}\II_{r_j}(x),
$$
and we get the expression similar to (\ref{mth}) with the only
difference that the outer sum is over all $\rho\in K$. However, if
$\rho\in K_t$ for $t>2r+2$, then by Theorem~\ref{th1},
$c(\rho)>r$, and hence $\II_r^\rho(x)\equiv 0$.
\end{proof}

Theorem~\ref{th4} provides a finite algorithm for finding
$\II_r(x)$ for any given $r\geq0$, since we have to consider all
involutions in $\II_{2r+2}$, and to perform certain routine
operations with all shapes found so far. Moreover, the amount of
search can be decreased substantially due to the following
proposition.

\begin{proposition}\label{ppp}
Let $\psi^0=21$, $\psi^1=3412$, $\psi^2=351624$, and
$$\psi^r=3\,\,5\,\,1\,\,7\,\,2\,\,9\,\,4\,\,11\ldots(2r+1)\,\,(2r-4)\,\,(2r+2)\,\,(2r-2)\,\,(2r)$$
for all $r\geq3$. Then the only kernel involution of capacity
$r\geq0$ and size $2r+2$ is $\psi^r$. Its contribution to
$\II_r(x)$ equals $\frac{x^{2r+2}}{(1-x)^r}\II_0^{r+2}(x)$.
\end{proposition}

This proposition is proved easily by induction, similarly to
Lemma~1. The feasible cells in the corresponding cell
decomposition is: $C_{ii}$ is free-cell if
$i=1,3,5,\ldots,2r+1,2r+2$, $C_{ii}$ is diagonal-decreasing-cell
if $i=2,4,\ldots,2r$, and all the other cells are infeasible
cells, hence the contribution to $\II_r(x)$ is as described. By
the above proposition, it suffices to search only involutions in
$\II_{2r+1}$. Below we present several explicit calculations.

\subsection{The case $r=0$}
Let us start from the case $r=0$. Observe that (\ref{mth}) remains
valid for $r=0$, provided the left hand side is replaced by
$\II_r(x)-1$; subtracting $1$ here accounts for the empty
involution. Also, by Subsection~2.2 we have only two shapes
$\rho^1=1$ and $\rho^2=21$ with $s(\rho^1)=1$, $c(\rho^1)=0$,
$f(\rho^1)=1$, $dd(\rho^1)=d(\rho^1)=0$, $s(\rho^2)=2$,
$c(\rho^2)=0$, $f(\rho^2)=2$, and $dd(\rho^2)=d(\rho^2)=0$.
Therefore, we get $$\II_0(x)-1=x\II_0(x)+x^2\II_0^2(x).$$

\begin{corollary}{\rm(see~\cite[Rem.~4.28]{GuibertThese})}\label{case0}
The generating function for the number of involutions which avoid
$3412$ is given by
    $$\II_0(x)=\frac{1-x-\sqrt{1-2x-3x^2}}{2x^2}.$$
\end{corollary}

\subsection{The case $r=1$} Let now $r=1$. The involutions in
$\II_4$ are exhibited only one kernel shape distinct from $\rho^1$
and $\rho^2$ which is $\rho^3=3412$, whose contribution equals
$\frac{x^4}{1-x}\II_0^3(x)$ (see Figure~4).
\begin{center}
\figd

{\sc Figure 4}. Kernel cell decomposition for $\pi\in
\II(1)\cup\II(21)\cup\II(3412)$.
\end{center}

Therefore, (\ref{mth}) amounts to
$\II_1(x)=x\II_1(x)+2x^2\II_0(x)\II_1(x)
+\frac{x^4}{1-x}\II_0^3(x)$, and we get the following result.

\begin{corollary}\label{case1}
The generating function for the number of involutions which
contain $3412$ exactly once is given by
$$\II_1(x)=-\frac{1-2x}{2x^2(1-x)}+\frac{1-2x-2x^2}{2x^2}\sqrt{1-2x-3x^2}^{\,-1}.$$
\end{corollary}

\subsection{The case $r=2$}
Let $r=2$. We have to check the kernel shapes of involutions in
$\II_6$. Exhaustive search adds one new shape to the previous
list; this is $\rho^4=351624$. Calculation of the parameters $s$,
$c$, $f$, $d$, $dd$ is straightforward, and we get

\begin{corollary}
The generating function for the number of involutions which
contain $3412$ exactly twice is given by
$$\II_2(x)=\frac{1-2x}{2x^2(1-x)}-\frac{1-6x+8x^2+8x^3-15x^4-2x^5+4x^6}{2x^2(1-x)^2}\sqrt{1-2x-3x^2}^{\,-3}.$$
\end{corollary}

\subsection{The cases $r=3,4,5,6,7$} Let $r=3,4,5,6$;
exhaustive search in $\II_{8}$, $\II_{10}$, $\II_{12}$, and
$\II_{14}$ reveals $2$, $5$, $12$, $25$, $48$, and $100$ new
kernel shapes, respectively, and we get

\begin{corollary}
Let $r=3,4,5,6,7$. Then the generating function for the number of
involutions which contain $3412$ exactly $r$ times is given by
$$\II_r(x)=\frac{1}{2x^2}F_r(x)+\frac{1}{2x^2}G_r(x)\sqrt{1-2x-3x^2}^{\,1-2r},$$
where

$\begin{array}{rl}
(1-x^2)F_3(x)&=-(1-2x)(1+x+x^2),\\
\\[-8pt]
(1-x^2)F_4(x)&=-1+3x+4x^2-8x^3-2x^4,\\
\\[-8pt]
(1-x^2)F_5(x)&=3-7x-7x^2+12x^3+6x^4,\\
\\[-8pt]
(1-x^2)^2F_6(x)&=-5+9x+21x^2-25x^3-34x^4+16x^5+24x^6-2x^7-2x^8,\\
\\[-8pt]
(1-x^2)^2F_7(x)&=7-11x-28x^2+20x^3+54x^4-2x^5-46x^6+2x^8,
\end{array}$

and

$\begin{array}{rl}
(1-x)^2G_3(x)&=1-8x+18x^2+x^2-29x^4-12x^5+14x^6+41x^7+2x^8-18x^9,\\
\\[-8pt]
(1-x)^4G_4(x)&=1-14x+71x^2-124x^3-166x^4+874x^5-624x^6-1332x^7+1909x^8\\
&+426x^9-1585x^{10}+292x^{11}+400x^{12}-126x^{13},\\
\\[-8pt]
(1-x)^4G_5(x)&=-3+46x-267x^2+627x^3+134x^4-3321x^5+3954x^6+5214x^7\\
             &-11775x^8-2186x^9+14525x^{10}-1701x^{11}-8824x^{12}+1537x^{13}\\
             &+2594x^{14}-216x^{15}-324x^{16},\\
\\[-8pt]
(1-x)^6G_6(x)&=5-94x+712x^2-2582x^3+3124x^4+8364x^5-31620x^6+15464x^7\\
             &+77508x^8-107098x^9-76814x^{10}+214160x^{11}+5782x^{12}-231050x^{13}\\
             &+62700x^{14}+146176x^{15}-65653x^{16}-50328x^{17}+29646x^{18}\\
             &+6462x^{19}-5346x^{20}+486x^{21},\\
\\[-8pt]
(1-x)^6G_7(x)&=-7+144x-1210x^2+5020x^3-8206x^4-12180x^5+69464x^6\\
             &-54210x^7-181468x^8+315366x^9+239852x^{10}-779338x^{11}\\
             &-124766x^{12}+1226006x^{13}-168810x^{14}-1272344x^{15}+418555x^{16}\\
             &+813368x^{17}-373802x^{18}-279554x^{19}+153648x^{20}+37188x^{21}\\
             &-23166x^{22}+486x^{23}.
\end{array}$
\end{corollary}
\section{Further results}
As an easy consequence of Theorem~\ref{th4} and
Corollary~\ref{case0} we get the following result.

\begin{corollary}\label{caa} Let $r\geq0$, then $\II_r(x)$ is a rational
function in the variables $x$ and $\sqrt{1-2x-3x^2}$.
\end{corollary}

Another direction would be to match the approach of this note with
the previous results on even (odd) permutations which contain
$132$ exactly $r$ times (see~\cite{Meven}). We denote the
generating function for the number of even (resp. odd) involutions
in $\II_n$ which contain $r$ occurrences of $3412$ by $E_r(x)$
(resp. $O_r(x)$) (see Table~\ref{teven}).
\begin{center}
\begin{table}[h]
\begin{tabular}{|c|c|c|c|c|c|c|c|c|c|c|c|c|c|}\hline
  $r\backslash n$ & $0$ & $1$ & $2$ & $3$ & $4$ & $5$ & $6$ & $7$ & $8$ & $9$ & $10$ & $11$ &
  $12$\\\hline\hline
  $0$ & $1$ & $1$ & $1$ & $2$ & $3$ & $11$ & $31$ & $71$ & $155$ & $379$ & $1051$ & $2971$ & $8053$\\ \hline
  $1$ & $0$ & $0$ & $0$ & $0$ & $1$ & $5$  & $14$ & $30$  & $82$ & $320$ & $1213$ & $3895$ & $11141$\\ \hline
  $2$ & $0$ & $0$ & $0$ & $0$ & $0$ & $0$  & $0$ & $0$  & $11$ & $95$ & $439$ & $1463$ & $4407$\\ \hline
  $3$ & $0$ & $0$ & $0$ & $0$ & $0$ & $0$  & $1$ & $4$  & $11$ & $29$ & $104$ & $396$ & $1486$\\ \hline
  $4$ & $0$ & $0$ & $0$ & $0$ & $0$ & $0$  & $0$ & $0$  & $14$ & $108$ & $321$ & $1612$ & $4782$\\ \hline
  $5$ & $0$ & $0$ & $0$ & $0$ & $0$ & $0$  & $0$ & $0$  & $6$ & $60$ & $275$ & $878$ & $2247$\\ \hline
  $6$ & $0$ & $0$ & $0$ & $0$ & $0$ & $0$  & $0$ & $0$  & $1$ & $21$ & $122$ & $446$ & $1504$\\ \hline
\end{tabular}
\caption{Number of even involutions in $\II_n$ which contain
$3412$ exactly $r$ times where $0\leq n\leq 12$ and $0\leq r\leq
6$.}\label{teven}
\end{table}
\end{center}
We define $N_r(x)=E_r(x)-O_r(x)$ for any $r$, that is,
$$N_r(x)=\sum_{n\geq0}\,\,\sum_{\pi\in\II_n\ \mbox{\small contains $3412$ exactly $r$ times}} (-1)^{21(\pi)}x^{n},$$
where $21(\pi)$ is the number of occurrences of the pattern $21$
in $\pi$. Using the arguments in the proof of Theorem~\ref{th4} we
get the following result.

\begin{theorem}\label{theven} For any $r\geq0$,
\begin{equation}
N_r(x)-\delta_{r,0}=\sum_{\rho\in
K_{2r+2}}\left[\frac{(-1)^{21(\rho)}x^{s(\rho)}(1+x)^{dd(\rho)}}{(1+x^2)^{d(\rho)+dd(\rho)}}
\sum_{r_1+\dots+r_{f(\rho)}=r-c(\rho)}
\;\prod_{j=1}^{f(\rho)}N_{r_j}(x)\right],\label{eventh}
\end{equation}
where $r_j\geq0$ for $1\leq j\leq f(\rho)$, $\delta_{0,0}=1$, and
$\delta_{r,0}=0$ for all $r\geq1$.
\end{theorem}

As an easy consequence of Theorem~\ref{theven} we get the
following result.

\begin{corollary}\label{cbb} Let $r\geq0$, then $N_r(x)$ is a rational
function in the variables $x$ and $\sqrt{1-2x+5x^2}$.
\end{corollary}
Clearly,
\begin{equation}
E_r(x)=\frac{1}{2}(\II_r(x)+N_r(x))\mbox{ and }
O_r(x)=\frac{1}{2}(\II_r(x)-N_r(x)), \label{eqeo}\end{equation}
for all $r\geq0$. Hence, Corollary~\ref{caa} and
Corollary~\ref{cbb} give the following result.
\begin{corollary} For any  $r\geq0$, the generating functions $E_r(x)$ and $O_r(x)$ are rational
functions in the variables $x$, $\sqrt{1-2x-3x^2}$, and
$\sqrt{1-2x+5x^2}$.
\end{corollary}
Hence, Theorem~\ref{theven} and Theorem~\ref{th4} provide a finite
algorithm for finding $\II_r(x)$, $N_r(x)$, $E_r(x)$, and $O_r(x)$
for any given $r\geq0$, since we have to consider all involutions
in $\II_{2r+2}$, and to perform certain routine operations with
all shapes found so far. Moreover, the amount of search can be
decreased substantially due to the following proposition which
holds as easily consequence of Proposition~\ref{ppp}.
\begin{proposition}
Let $\psi^0=21$, $\psi^1=3412$, $\psi^2=351624$, and
$$\psi^r=3\,\,5\,\,1\,\,7\,\,2\,\,9\,\,4\,\,11\ldots(2r+1)\,\,(2r-4)\,\,(2r+2)\,\,(2r-2)\,\,(2r)$$
for all $r\geq3$. Then the only kernel involution of capacity
$r\geq0$ and size $2r+2$ is $\psi^r$. Its contribution to $N_r(x)$
equals $(-1)^{r+1}\frac{x^{2r+2}(1+x)^r}{(1+x^2)^r}N_0^{r+2}(x)$.
\end{proposition}
By the above proposition, it suffices to search only involutions
in $\II_{2r+1}$. Similarly as our calculations for $\II_r(x)$
where $r=0,1,2,3,4,5,6,7$ with using Theorem~\ref{theven} we get
the following result.
\begin{corollary}\label{cev}
Let $r=0,1,2,3,4,5,6,7$. Then the generating function $N_r(x)$ is
given by
$$N_r(x)=\frac{1}{2x^2}P_r(x)+\frac{1}{2x^2}Q_r(x)\sqrt{1-2x+5x^2}^{\,1-2r},$$
where $$\begin{array}{rl}
P_0(x)&=x-1,\\
\\[-9pt]
(x^2+1)P_1(x)&=(x+1)(2x^2-2x+1),\\
\\[-9pt]
(x^2+1)^2P_2(x)&=(x-1)(2x^2-2x+1)(1+x)^2,\\
\\[-9pt]
(x^2+1)^3P_3(x)&=(2x^2-2x+1)(x^6+x^5+3x^4-2x^3-x^2+x+1),\\
\\[-9pt]
(x^2+1)^4P_4(x)&=(1-x^2)(2x^8-2x^7+10x^6+x^5+15x^4-12x^3+12x^2-3x+1),\\
\\[-9pt]
(x^2+1)^5P_5(x)&=-4x^{13}-14x^{11}+5x^{10}-33x^9+29x^8-16x^7+34x^6-42x^5+14x^4+6x^3\\
&-15x^2+7x-3,\\
\\[-9pt]
(x^2+1)^6P_6(x)&=-6x^{16}+6x^{15}-44x^{14}+58x^{13}-128x^{12}+163x^{11}-195x^{10}+271x^9-221x^8\\
&+188x^7-170x^6+160x^5-82x^4+27x^3+9x^2-9x+5,\\
\\[-9pt]
(x^2+1)^7P_7(x)&=-14x^{18}+12x^{17}-80x^{16}+94x^{15}-176x^{14}+212x^{13}-188x^{12}+247x^{11}\\
&-157x^{10}+35x^9-51x^8+28x^7+62x^6-142x^5+102x^4-49x^3-3x^2+11x-7,
\end{array}$$
and
$$\begin{array}{rl}
Q_0(x)&=1,\\
\\[-9pt]
(x^2+1)Q_1(x)&=(x^2-1)(4x^2-2x+1),\\
\\[-9pt]
(x^2+1)^2Q_2(x)&=(22x^6-58x^5+69x^4-48x^3+22x^2-6x+1)(1+x)^2,\\
\\[-9pt]
(x^2+1)^3Q_3(x)&=(x-1)\bigl(100x^{12}-18x^{11}+323x^{10}-507x^9+491x^8-182x^7+52x^6-14x^5\\
&+46x^4-34x^3+19x^2-5x+1\bigr),\\
\\[-9pt]
(x^2+1)^4Q_4(x)&=-(1+x)\bigl(650x^{16}-1880x^{15}+5992x^{14}-9143x^{13}+13671x^{12}-19666x^{11}\\
&+26606x^{10}-28683x^9+24771x^8-16778x^7+9158x^6-3969x^5+1385x^4-374x^3\\
&+78x^2-11x+1\bigr),\\
\\[-9pt]
(x^2+1)^5Q_5(x)&=(1-x)\bigl(5000x^{21}-4650x^{20}+24624x^{19}-25585x^{18}+76987x^{17}-95269x^{16}\\
&+127936x^{15}-140244x^{14}+169896x^{13}-159580x^{12}+119898x^{11}-51878x^{10}\\
&-84x^9+26302x^8-26778x^7+17822x^6-8604x^5+3270x^4-940x^3+209x^2\\
&-31x+3\bigr),
\end{array}$$
$$\begin{array}{rl}
(x^2+1)^6Q_6(x)&=-43750x^{27}+133750x^{26}-542500x^{25}+1329674x^{24}-2984612x^{23}\\
                &+5378699x^{22}-8590394x^{21}+12236909x^{20}-15828644x^{19}+18229621x^{18}\\
                &-19177696x^{17}+18659837x^{16}-17024788x^{15}+14266232x^{14}-10700428x^{13}\\
                &+6908636x^{12}-3700402x^{11}+1527142x^{10}-395516x^9-27686x^8+101584x^7\\
               &-68679x^6+30486x^5-10181x^4+2580x^3-493x^2+64x-5,\\
\\[-8pt]
(x^2+1)^7Q_7(x)&=-481250x^{31}+1658750x^{30}-5844500x^{29}+14332172x^{28}-29824134x^{27}\\
                &+52203592x^{26}-78380980x^{25}+104774831x^{24}-124983968x^{23}+132048678x^{22}\\
                &-122776812x^{21}+101431782x^{20}-72478438x^{19}+39230434x^{18}-4374004x^{17}\\
                &-25236483x^{16}+42491126x^{15}-44337242x^{14}+34831062x^{13}-21298364x^{12}\\
               &+9941638x^{11}-3111220x^{10}+199166x^9+519349x^8-425180x^7+212566x^6\\
               &-78950x^5+22882x^4-5138x^3+874x^2-102x+7.
\end{array}$$
\end{corollary}

For example, by Corollary~\ref{case0} and Corollary~\ref{cev} for
$r=0$ we have that the generating function for the number even
involutions which avoid $3412$ is given by
   $$E_0(x)=\frac{\sqrt{1-2x+5x^2}-\sqrt{1-2x-3x^2}}{4x^2},$$
and the generating function for the number of odd involution which
avoid $3412$ is given by
   $$O_0(x)=\frac{2-2x-\sqrt{1-2x+5x^2}-\sqrt{1-2x-3x^2}}{4x^2}.$$


\begin{thebibliography}{99}
\bibitem{AlFr00}
N.~Alon and E.~Friedgut, On the number of permutations avoiding a
given pattern, {\em J. Combin. Theory Ser. A} 89(1):133--140,
2000.

\bibitem{At99}
M.~D. Atkinson, Restricted permutations, {\em Discrete Math.}
195(1-3):27--38, 1999.

\bibitem{Bo97b}
M.~B{\'o}na, Exact enumeration of $1342$-avoiding permutations: a
close link with labeled trees and planar maps, {\em J. Combin.
Theory Ser. A} 80(2):257--272, 1997.

\bibitem{Bo97c}
M.~B{\'o}na, The number of permutations with exactly $r$
$132$-subsequences is ${P}$-recursive in the size!, {\em Adv. in
Appl. Math.} 18(4):510--522, 1997.

\bibitem{Bo97a}
M.~B{\'o}na, Permutations avoiding certain patterns: the case of
length $4$ and some generalizations, {\em Discrete Math.}
175(1-3):55--67, 1997.

\bibitem{Bo98}
M.~B{\'o}na, Permutations with one or two $132$-subsequences, {\em
Discrete Math.} 181(1-3):267--274, 1998.

\bibitem{ChWe99}
T.~Chow and J.~West, Forbidden subsequences and {C}hebyshev
polynomials, {\em Discrete Math.} 204(1-3):119--128, 1999.

\bibitem{Egge}
E.S.~Egge, Restricted $3412$-avoiding involutions: Continued
fractions, Chebyshev polynomials and enumerations, {\em Adv. in
Appl. Math.}, to appear.

\bibitem{EM}
E.S.~Egge and T.~Mansour, Involutions restricted by $3412$,
Continued fractions, and Chebyshev polynomials, preprint.

\bibitem{Fulmek}
M.~Fulmek, Enumeration of permutations containing a presribed
number of occurrences of a pattern of length three, {\em Adv.
Appl. Math.} 30:607--632, 2003.

\bibitem{Gessel}
I.M.~Gessel, Symmetric functions and P-recursiveness, {\em J.
Comb. Th. A} 53:257--285, 1990.

\bibitem{GireThese}
S.~Gire, Arbres, permutations \`a motifs exclus et cartes
planaires~: quelques probl\`emes algorithmiques et combinatoires,
{\em PHD-thesis, University Bordeaux~1, France} (1993).

\bibitem{GouyouBYoung}
D.~Gouyou--Beauchamps, Standard Young tableaux of height $4$ and
$5$, {\em European J. Combin.} 10:69--82, 1989.

\bibitem{GuibertThese}
O.~Guibert, Combinatoire des permutations \`a motifs exclus en
liaison avec mots, cartes planaires et tableaux de Young, {\em
PHD-thesis, University Bordeaux~1, France} (1995).

\bibitem{GPP}
O.~Guibert, E.~Pergola and R.~Pinzani, Vexillary involutions are
enumerated by Motzkin numbers, {\em Annals of Comb.} 5:153--174,
2001.

\bibitem{Jag}
A.D.~Jaggard, Prefix Exchanging and pattern avoidance by
involutions, prepint math.CO/0306002.

\bibitem{Kn}
D.E.~Knuth, The Art of Computer Programming, 2nd ed. Addison
Wesley, Reading, MA (1973).

\bibitem{Ma00}
T.~Mansour, Permutations containing and avoiding certain patterns,
In {\em Formal power series and algebraic combinatorics (Moscow,
  2000)}, 704--708. Springer, Berlin, 2000.

\bibitem{Meven}
T.~Mansour, Counting occurrences of $132$ in an even permutation,
{\em International Journal of Mathematics and Mathematical
Sciences}, to appear, preprint math.CO/0211205.

\bibitem{MaVa01}
T.~Mansour and A.~Vainshtein, Restricted 132-avoiding
permutations, {\em Adv. Appl. Math.} 26:258--269, 2001.

\bibitem{MaVa02}
T.~Mansour and A.~Vainshtein, Counting occurrences of $132$ in a
permutation, {\em Adv. Appl. Math.} 28(2):185--195, 2002.

\bibitem{No96}
J.~Noonan, The number of permutations containing exactly one
increasing subsequence of length three, {\em Discrete Math.}
152(1-3):307--313, 1996.

\bibitem{NoZe96}
J.~Noonan and D.~Zeilberger, The enumeration of permutations with
a prescribed number of ``forbidden'' patterns, {\em Adv. in Appl.
Math.} 17(4):381--407, 1996.

\bibitem{Regev}
A.~Regev, Asymptotic values for degrees associated with strips of
Young diagrams, {\em Adv. Math.} 41:115-136, 1981.

\bibitem{Ro99}
A.~Robertson, Permutations containing and avoiding $123$ and $132$
patterns, {\em Discrete Math. Theor. Comput. Sci.} 3(4):151--154
(electronic), 1999.

\bibitem{SiSc85}
R.~Simion and F.~W. Schmidt, Restricted permutations, {\em
European J. Combin.} 6(4):383--406, 1985.

\bibitem{St94}
Z.~Stankova, Forbidden subsequences, {\em Discrete Math.}
132(1-3):291--316, 1994.

\bibitem{St96}
Z.~Stankova, Classification of forbidden subsequences of length
$4$, {\em European J. Combin.} 17(5):501--517, 1996.

\bibitem{We95}
J.~West, Generating trees and the {C}atalan and {S}chr\"oder
numbers, {\em Discrete Math.} 146(1-3):247--262, 1995.
\end{thebibliography}
\end{document}